# Solving Polynomial Systems by Penetrating Gradient Algorithm Applying Deepest Descent Strategy

Nikica Hlupić [1]   and   Ivo Beroš [2]

[1] University of Zagreb, Faculty of Electrical Engineering and Computing, Unska 3,
Zagreb 10000, Croatia; e-mail: nikica.hlupic@fer.hr
[2] VERN, University of Applied Sciences, Trg bana Josipa Jelačića 3, Zagreb 10000, Croatia;
e-mail: ivo.beros@vern.hr

**Abstract.** An algorithm and associated strategy for solving polynomial systems within the optimization framework is presented. The algorithm and strategy are named, respectively, *the penetrating gradient algorithm* and *the deepest descent strategy*. The most prominent feature of *penetrating gradient algorithm*, after which it was named, is its ability to "see and penetrate through" the obstacles in error space along the line of search direction and to jump to the global minimizer in a single step. The ability to find the deepest point in an arbitrary direction, no matter how distant the point is and regardless of the relief of error space between the current and the best point, motivates movements in directions in which cost function can be maximally reduced, rather than in directions that seem to be the best locally (like, for instance, *the steepest descent*, i.e., negative gradient direction). Therefore, the strategy is named *the deepest descent*, in contrast but alluding to *the steepest descent*. *Penetrating gradient algorithm* is derived and its properties are proven mathematically, while features of *the deepest descent* strategy are shown by comparative simulations. Extensive benchmark tests confirm that the proposed algorithm and strategy jointly form an effective solver of polynomial systems. In addition, further theoretical considerations in Section 5 about solving linear systems by the proposed method reveal a surprising and interesting relation of proposed and Gauss-Seidel method.

## 1.    Introduction

Solving systems of nonlinear equations is of high practical importance because they frequently appear in many fields, like mechanical and civil engineering or robotics and control theory, to mention few of them. Typically, it is not known whether the exact solution exists at all, so usual approach is to set the problem like the problem of minimization of residual sum of squares (*rss*) in the system and to apply an optimization method to find the best possible solution according to least squares (LS) principle. There are several methods [1] that can be tried to solve the problem but, generally, if equations are differentiable, suitable methods are so-called *line search* methods like gradient, Newton's (Levenberg-Marquardt variant), Quasi-Newton and possibly conjugate gradient [1],[2], while with non-differentiable equations we use derivative-free methods such as Nelder-Mead [1],[2],[3]. Finally, if feasible set (possible values of unknowns) is closed [2] or we, at least, can assign probabilities (goodness) to candidate points, various heuristics (simulated annealing, particle swarm, etc.; also called *stochastic* algorithms) and genetic algorithms can be effective as well [2],[3]. Because the underlying theory is best developed and convergence guaranteed only for derivative-based methods, they are usually the most desirable ones. Unfortunately, their actual implementation is not as simple as in theory because in a general unconstrained optimization problem we are confronted with two issues: first, we have to choose search direction and second, we have to determine how far to move (i.e., step). Answering this second question requires line search, which is still a bottle-neck of optimization algorithms and a research challenge [4],[5],[6],[7],[28]. However, it turns out that with polynomial systems we have an attractive alternative.

There are several approaches to solving polynomial systems, ranging from "classical" optimization [8] to branch and prune approach [9],[10], symbolic algorithms [11] and homotopy continuation [12],[13] but none of them is absolutely superior for all purposes. The method presented in this paper fits within the optimization framework but is distinctive by its underlying computation that avoids line search and serves as the engine for a specific strategy that we call *the deepest descent* strategy. The computational engine is named *penetrating gradient algorithm* after its most prominent feature, i.e., its ability to see and penetrate through the obstacles in the error space. The penetration is achieved by a remarkable property of the algorithm that it searches globally and finds the exact global minimizer of the objective function along the whole line of search direction (the whole $\mathfrak{R}$). This means that if the best point on the line of search direction is "behind" some obstacle, the algorithm sees it anyway and penetrates through the obstacle to move directly to the minimizer, no matter how distant it is from the current position.

The fact that we are able to find the deepest point of the error space along an arbitrary search direction enables and motivates a "relaxed" optimization strategy as simple as the following: determine suitable (promising) directions, find the deepest point on each of these lines and move to the deepest point found. In contrast, but alluding to, the steepest descent method, a natural title for this strategy is *the deepest descent*.

The ability to "jump" all over the whole $\Re$ usually results in faster convergence than achievable by other methods, provided that we approach the optimum from "the same side". This can be illustrated by extreme (fortuitous) case that we accidentally (or not) direct the search right to the solution, i.e., the solution lays on the line of search direction. The solution can be very far away and behind obstacles, and as such invisible to other search methods, but penetrating gradient algorithm will find it in a single step. In such situation, the deepest descent method would solve the system instantaneously. Moreover, the ability to see and penetrate through the obstacles makes the proposed method not only very fast, but also very successful in searching for a global minimizer, i.e., solution. This is a direct consequence of the penetrating feature of the algorithm because it is less likely that the deepest point will remain hidden after "piercing" the error space by many "rays" during the deepest descent optimization. All these facts will become apparent by mathematical derivation of penetrating gradient algorithm in the next section and comparative tests of the proposed method by simulations in the later sections of the paper.

## 2. The Penetrating Gradient Algorithm

Consider a system of $m$ nonlinear equations in $n$ unknowns $\mathbf{x} = [x_1\ x_2\ \ldots\ x_n]^T \in R^n$, where we denote all equations as $F_i(\mathbf{x})$, $F_i : R^n \rightarrow R$. With notation $\mathbf{F}: R^n \rightarrow R^m$, $\mathbf{F}(\mathbf{x}) = [F_1(\mathbf{x})\ F_2(\mathbf{x})\ \ldots\ F_m(\mathbf{x})]^T$ and $\mathbf{b} = [b_1\ b_2\ \ldots\ b_m]^T \in R^m$, the system can be written as

$$\mathbf{F}(\mathbf{x}) = \mathbf{b}. \qquad (2.1)$$

Residuals in the system are $r_i = b_i - F_i(\mathbf{x})$, i.e., vector of residuals is $\mathbf{r} = [r_1\ r_2\ \ldots\ r_m]^T \in R^m$. Commonly, even if the system has no exact solution, we wish to find the best possible solution according to least squares, so the goal is to minimize residual sum of squares

$$rss(\mathbf{x}) = \sum [b_i - F_i(\mathbf{x})]^2 = \|\mathbf{b} - \mathbf{F}(\mathbf{x})\|^2 = \|\mathbf{r}\|^2 \geq 0 \quad, \qquad (2.2)$$

where at the exact solution is found when $rss = 0$. To this aim we rely on the first order necessary condition [2] for a minimizer and we search for a solution to the equation

$$\nabla rss(\mathbf{x}) = \mathbf{0}. \qquad (2.3)$$

However, equation (2.3) is rarely solvable algebraically so we resort to optimization. The optimization approach is, in general, to minimize $rss(\mathbf{x})$ applying an iterative procedure in which we take the current point (estimate of solution) $\mathbf{x}$, decide about a search direction $\mathbf{d} \in R^n$ and then search for the step $\alpha \in R$ such that new point given by $\mathbf{x} + \alpha \cdot \mathbf{d}$ yields $rss(\mathbf{x} + \alpha \cdot \mathbf{d})$ as low as possible. Symbolically, the requirement is

$$\phi(\alpha) = rss(\mathbf{x} + \alpha \cdot \mathbf{d}) = \|\mathbf{b} - \mathbf{F}(\mathbf{x} + \alpha \cdot \mathbf{d})\|^2 \rightarrow min \qquad (2.4)$$

and we need $\alpha^*$ such that

$$\alpha^* = \arg\min_{\alpha \in R} \phi(\alpha) = \arg\min_{\alpha \in R} \|\mathbf{b} - \mathbf{F}(\mathbf{x} + \alpha \cdot \mathbf{d})\|^2 = \arg\min_{\alpha \in R} \sum_{i=1}^{m} [b_i - F_i(\mathbf{x} + \alpha \cdot \mathbf{d})]^2 . \qquad (2.5)$$

Search for $\alpha^*$ in (2.5) is one-dimensional search (i.e., line search) and in general it is not a simple task. Therefore, here we depart from "classical" approach and avoid line search. The alternative follows from the fact that $\mathbf{F}(\mathbf{x})$ is a polynomial system [12], i.e., all $F_i$ are polynomials, meaning that all $F_i(\mathbf{x})$ in (2.1) are linear combinations of terms of the form

$$a \cdot x_1^{p1} x_2^{p2} \ldots x_n^{pn} \quad, \qquad (2.6)$$

where $a$ is a real number and $p_i$ are nonnegative integers (for $p_i = 0$, variable $x_i$ does not figure in a term). Let us, for convenience, write the system (2.1) with all $F_i$ being polynomials as

$$P(\mathbf{x}) = \mathbf{c} \quad . \qquad (2.7)$$

With all $F_i$ being polynomials and notation (2.7), the problem (2.5) of line search becomes

$$\alpha^* = \arg\min_{\alpha \in R} \phi(\alpha) = \arg\min_{\alpha \in R} \|\mathbf{c} - P(\mathbf{x} + \alpha \cdot \mathbf{d})\|^2 = \arg\min_{\alpha \in R} \sum_{i=1}^{m} [c_i - P_i(\mathbf{x} + \alpha \cdot \mathbf{d})]^2 \quad . \qquad (2.8)$$

Now we make a key observation. Namely, because the current point $\mathbf{x}$ and search direction $\mathbf{d}$ are known, all arguments (variables) of polynomials $P_i(\mathbf{x}+\alpha \cdot \mathbf{d})$ in (2.8) are just first degree polynomials in $\alpha$. This is obvious if we recall that in (2.8) we have polynomials $P_i(x_1+\alpha \cdot d_1, x_2+\alpha \cdot d_2, \ldots, x_n+\alpha \cdot d_n)$, where $\alpha$ is the only unknown.

Furthermore, since all $P_i$ are constituted of terms of the form (2.6), that is, products of arguments (possibly raised to nonnegative integer powers), in (2.8) all $P_i$ themselves become polynomials in $\alpha$. For example, if $i^{th}$ equation in the system (2.7) is

$$P_i: \quad x_1 x_2 + x_1^2 x_2 = c_i \quad ,$$

in (2.8) it will contribute as

$$[c_i - (x_1 + \alpha \cdot d_1)(x_2 + \alpha \cdot d_2) - (x_1 + \alpha \cdot d_1)^2 \cdot (x_2 + \alpha \cdot d_2)]^2 \quad . \tag{2.9}$$

Because $\mathbf{x} = [x_1 \ x_2]^T$ and $\mathbf{d} = [d_1 \ d_2]^T$ are known, the term (2.9) is simply a polynomial in $\alpha$ (sixth degree polynomial). Hence, in (2.8) we have summation of polynomials in $\alpha$ and it follows that the whole sum (2.8) is nothing but a polynomial (in $\alpha$), which we call "summary polynomial". Since polynomials are easily differentiable and their roots can be computed by well defined and efficient routines [26], it follows that actually there is no need for line search to minimize (2.8). Instead, minimization of (2.8) with respect to $\alpha$ can be achieved algebraically by solving

$$\phi'(\alpha) = d\phi(\alpha) / d\alpha = 0, \tag{2.10}$$

where $\phi(\alpha)$ is from (2.8), which leads to

$$\phi'(\alpha) = \sum_{i=1}^{m} [c_i - P_i(\mathbf{x} + \alpha \cdot \mathbf{d})] \frac{dP_i(\mathbf{x} + \alpha \cdot \mathbf{d})}{d\alpha} = 0 \quad . \tag{2.11}$$

Terms $(c_i - P_i)$ in (2.11) are polynomials in $\alpha$ of the same degree as $P_i$ (for the definition of the degree see, for example, [12]). Let us denote the degree of $P_i$ as $\beta_i$. The derivative terms $dP_i/d\alpha$ are polynomials of the degree one less than the degree of $P_i$ so their degrees are $(\beta_i - 1)$. Hence, the sum in (2.11) consists of products of polynomials and is again just a $\beta + (\beta - 1) = 2\beta - 1$ polynomial, where $\beta = max\{\beta_1, \ldots, \beta_m\}$. This polynomial certainly has negative leading coefficient because in (2.11) we multiply $-P_i \cdot dP_i/d\alpha$. It follows that the best step $\alpha^*$ we search for is simply one of the real roots of polynomial (2.11). Clearly, the degree $2\beta - 1$ of $\phi'(\alpha)$ is an odd number, so there will be at least one real root, which conforms to intuitive expectation that a real number $\alpha$ minimizing $rss$ given by (2.11) certainly exists.

By solving (2.11) we obtain $\alpha^*$ directly, i.e., without line search. Moreover, since we solve $\phi'(\alpha) = 0$, the result will be the best possible $\alpha$ in the whole R, so the effect is as if we searched interval $(-\infty, \infty)$, which is unachievable by any line search routine. Unfortunately, the solution to $\phi'(\alpha) = 0$ defines stationary points, not necessarily minimizers, of $\phi(\alpha)$ function. The resulting equation (2.11) is $2\beta - 1$ degree polynomial and there will be $2\beta - 1$ roots but at present we do not know how to identify the right one in advance. Therefore, we compute all roots, determine $rss$ for each real root and eventually set $\alpha^*$ to be the root yielding the lowest $rss$. In practice, the degree of summary polynomial can easily reach few tens so finding its minimizer (i.e., solution to $\phi'(\alpha) = 0$) is numerically the most demanding computation in the proposed algorithm. The roots-finding procedure itself is not the subject of this paper and the interested reader should consult the literature on numerical analysis [14] (Matlab®, for example, computes the roots as eigenvalues of the companion matrix). The determination of the right root in advance is certainly a challenge for future research and improvements. Now we are ready to formulate penetrating gradient algorithm for polynomial systems.

**Algorithm 1:** Penetrating gradient algorithm for polynomial systems
1. Determine (in any way; irrelevant at this point) search direction $\mathbf{d}$.
2. Given the current point $\mathbf{x}$ and search direction $\mathbf{d}$, substitute $(\mathbf{x} + \alpha \cdot \mathbf{d})$ for $\mathbf{x}$ in (2.1), that is, in (2.7). This will render the system if polynomials in $n$ variables $\mathbf{x} = [x_1 \ x_2 \ \ldots \ x_n]^T \in R^n$ into a system of polynomials in a single unknown $\alpha$.
3. Construct equation (2.11). It will be a polynomial in $\alpha$.
4. Compute the roots of polynomial (2.11) constructed in the step 3.
5. Among the roots obtained in the step 4, select a real one that yields the lowest $rss$ (lowest sum in (2.8)). The selected solution $\alpha^*$ is the best step $\alpha \in (-\infty, \infty)$ for movement along the line of search direction $\mathbf{d}$.
6. Compute new point $\mathbf{x}_{new} = \mathbf{x} + \alpha^* \cdot \mathbf{d}$. This $\mathbf{x}_{new}$ minimizes the $rss$ of the original system (2.7) along the line of search direction $\mathbf{d}$.

We conclude development by formally establishing the property of penetrating gradient algorithm that the point it yields is the best possible point, i.e., the global minimizer in a given search direction.

**Proposition 1:** *Let f(**x**) be the objective function, in this case residual sum of squares (rss). Starting from any point **x**$^{(k)}$ and searching in any direction, the point **x**$_{pg}^{(k+1)}$ obtained by penetrating gradient algorithm will be the global minimizer of f(**x**) in a given search direction. That is, f(**x**$_{pg}^{(k+1)}$) ≤ f(**x**) for all **x** along the line of search direction.* □

*Proof*: Let **x**$^{(k)}$ be the current point and **d** an arbitrary search direction. By step 6 of Algorithm 1, point **x**$_{pg}^{(k+1)}$ obtained by penetrating gradient algorithm is point **x**$_{pg}^{(k+1)}$ = **x**$_{new}$ = **x**$^{(k)}$ + $\alpha$*·**d**, where $\alpha$* is the best step among all candidate steps (by step 5 of Algorithm 1), i.e., among all $\alpha$ that (by step 4 of Algorithm 1) satisfy the first order necessary condition (2.10) for a minimizer of f(**x**$^{(k)}$+$\alpha$·**d**). Thus, f(**x**$_{pg}^{(k+1)}$) = f(**x**$^{(k)}$+$\alpha$*·**d**) ≤ f(**x**$^{(k)}$+$\alpha$·**d**) = f(**x**$^{(k+1)}$) for all **x**$^{(k+1)}$ = **x**$^{(k)}$+$\alpha$·**d**, $\alpha \in (-\infty, \infty)$. ∎

In regard to the name of the algorithm, it is motivated by its properties. Line search has been replaced by computation of LS solution for a minimizer of *rss* of the system of polynomials (2.7) along the line of search direction, i.e., by computation of "the best" real root of polynomial (2.11). The fact that LS solution is the best possible solution in the whole set $\Re$ means that we, effectively, can see and move from $-\infty$ to $\infty$ as if we were able to penetrate through the obstacles and move right to the global minimizer (along search direction). This is strong motivation for the title *penetrating*.

Ability to work with arbitrary search directions makes this algorithm a powerful computational engine for minimization of *rss*, but success in finding solution to system (2.7) depends on recognition of appropriate search directions. We may have various indications of promising directions but whichever directions we try, the least expected from any optimization method is that it exhibits global convergence (i.e. approaches a stationary point of the objective function [1],[2]). The easiest way to ensure global convergence is to "inherit" it from gradient method, which promotes the direction of (negative) gradient to an inevitable search direction. This fact justifies complete name *penetrating gradient algorithm*, despite the fact that in actual performance of the algorithm the gradient direction is quite a rare choice, as explained in the next section.

### 3. The deepest descent strategy

In addition to the direction of gradient, which is an inevitable one as already explained, in our work we also decided to regularly try directions of all coordinate axes (similarly to coordinate descent method [30]) as well as the Newton's direction (when number of equations in (2.7) is equal to the number of unknowns). In this way we pierce the space symmetrically and in most "relevant" directions, hoping that thereby we reduce the possibility to miss the global optimum, i.e., solution. Having obtained the best points in several directions, we follow greedy strategy [18] and accept movement along direction in which the achievable reduction of *rss* is maximal. Consequently, in each iteration we move to the deepest detected point and therefore we named this strategy *the deepest descent*.

**Algorithm 2:** The deepest descent method
1. Given the current point **x**, compute the gradient **g**(**x**) of the objective function (in our discussion *rss*). Set search direction to **d** = –**g**(**x**).
2. Apply penetrating gradient algorithm (Algorithm 1) with **d** (from step 1) as search direction to find the global minimizer **x**$_g$ of the objective function along the line of gradient.
3. Set **x**$_{best}$ = **x**$_g$ and rss$_{min}$ = *rss*(**x**$_g$) as referential values for further use.
4. `for` all search directions we wish to try (Newton's, coordinate axes, …)
    a) Set **d** = search direction.
    b) Apply penetrating gradient algorithm (Algorithm 1) with **d** (from step 4a) as search direction to find the global minimizer **x**$_d$ of the objective function along the line of search direction **d**.
    c) If *rss*(**x**$_d$) < rss$_{min}$
        set **x**$_{best}$ = **x**$_d$ and rss$_{min}$ = *rss*(**x**$_d$).
5. Move to point **x**$_{best}$, i.e., set **x** = **x**$_{best}$.
6. Check stopping criteria and if none is satisfied, continue from the beginning (step 1).

Though, strictly speaking, Algorithm 2 uses both the penetrating gradient algorithm and the deepest descent strategy, for convenience we named Algorithm 2 *The Deepest Descent Method*.

Now, based on the fact that one of search directions is necessarily the direction of (negative) gradient, we can formally establish the descent property of the proposed method.

**Proposition 2:** *Let $\mathbf{x}_{dd}$ be a point obtained by deepest descent method according to Algorithm 2. If $\{\mathbf{x}_{dd}^{(k)}\}_{k=0}^{\infty}$ is the deepest descent sequence for $f : \mathbb{R}^n \rightarrow \mathbb{R}$ and if $\nabla f(\mathbf{x}^{(k)}) \neq 0$, then $f(\mathbf{x}_{dd}^{(k+1)}) < f(\mathbf{x}_{dd}^{(k)})$. In words, the deepest descent method exhibits descent property. Furthermore, for any starting point $\mathbf{x}^{(k)}$, the deepest descent method will yield a point either better or equal to the point obtained by the steepest descent method.* □

*Proof*: To prove both stated properties of the deepest descent method, recall that descent property holds for the steepest descent algorithm [1],[2]. Hence, we know that whenever $\nabla f(\mathbf{x}^{(k)}) \neq 0$, the steepest descent algorithm yields a point $\mathbf{x}_{sd}^{(k+1)}$ for which $f(\mathbf{x}_{sd}^{(k+1)}) < f(\mathbf{x}^{(k)})$. Furthermore, by Proposition 1 and the fact that one of search directions in the deepest descent method is necessarily the direction of (negative) gradient (by steps 1 and 2 in Algorithm 2), we also know that for any current point $\mathbf{x}^{(k)}$ it will always be $f(\mathbf{x}_{dd}^{(k+1)}) \leq f(\mathbf{x}_{sd}^{(k+1)})$. Thus, when $\nabla f(\mathbf{x}^{(k)}) \neq 0$, we have $f(\mathbf{x}_{dd}^{(k+1)}) \leq f(\mathbf{x}_{sd}^{(k+1)}) < f(\mathbf{x}^{(k)})$, which completes the proof. ∎

Simulations presented in the section with comparative study show that *the deepest descent method* successfully solves systems with several equations in several unknowns and it usually takes no more than few hundred iterations for *rss* to be sufficiently reduced. Nevertheless, with larger systems the algorithm can run into the state in which it advances very slowly. Deeper insight into the program indicates that in such situations residuals in the majority of equations decrease, while in one or two equations residuals increase. This is a state from which the algorithm hardly gets out and in order to reduce probability of running into such situations, we introduce one heuristics in the method. In addition to comparing *rss* reductions obtained in candidate points (step 4c in Algorithm 2), we also watch particular residuals in the system and we accept the point that maximally reduces *rss* but under the condition that none of the residuals in candidate point exceeds the largest residual in a referential point (which is the point found in gradient direction). Since *rss* reduction is the primary aim and maximum residual (norm infinity) a secondary control condition, in what follows we denote the Algorithm 2 with incorporated heuristics as *rss1rmax2*.

## 4. Comparative study

Although penetrating gradient algorithm requires specific class of equations, i.e., polynomial systems, it follows classical optimization approach to solving systems of nonlinear equations, that is, it tries to minimize the residual sum of squares. The aim is to find a real solution to the system if it exists, and if there are many real solutions, any one is acceptable. Therefore, the proposed method can be directly compared only with classical optimization methods because methods specialized for polynomial systems, like homotopy continuation, represent completely different approach with different objectives.

The deepest descent method has been extensively tested and compared with several confirmed algorithms, specifically with Newton's method (Levenberg-Marquardt variant), conjugate gradient and Nelder-Mead algorithm. However, conjugate gradient and Nelder-Mead showed to be so much less effective with polynomial systems than Newton's and deepest descent method that their results are not reported here for clarity and conciseness.

This comparative study firstly examines the basic algorithm properties and then it demonstrates the performance of the proposed method on many benchmark examples of polynomial systems.

**Demonstration of penetrating property**

In order to visualize the algorithm's penetrating property, we consider a well-known two-dimensional problem, i.e., we minimize Rosenbrock function, known to be of relatively simple mathematical form, but still quite a demanding one for optimization [2]. Rosenbrock function is defined by

$$f(\mathbf{x}) = 100(x_2 - x_1^2)^2 + (1 - x_1)^2 \tag{4.1}$$

and it is suitable for penetrating gradient algorithm, which becomes apparent if we rewrite (4.1) as

$$f(\mathbf{x}) = 100 x_1^4 + 100 x_2^2 + x_1^2 - 200 x_1^2 x_2 - 2 x_1 + 1 \;. \tag{4.2}$$

Rosenbrock function has a global minimum in $\mathbf{x} = [1 \; 1]^T$ and to find it we drive gradient $\mathbf{g}(\mathbf{x}) = \nabla f(\mathbf{x})$ to zero. Obviously, equating the gradient

$$\mathbf{g}(\mathbf{x}) = \begin{bmatrix} 400 x_1^3 + 2 x_1 - 400 x_1 x_2 - 2 \\ 200 x_2 - 200 x_1^2 \end{bmatrix} \tag{4.3}$$

with zero forms a system of two polynomial equations in two unknowns $x_1$ and $x_2$, which is a system of the form (2.7) as required by penetrating gradient algorithm.

Figure 1 illustrates the penetrating property of the algorithm. We start from $\mathbf{x}^{(0)} = [0\ 0]^T$ and compare pure gradient method that performs line search by backtracking with penetrating gradient algorithm that searches only in the direction of the gradient. As apparent in the figure, starting point and optimum are on different sides of the hill. Blue line is trajectory of pure gradient method and we see that it approaches optimum very slowly, sliding along the surface defined by Rosenbrock function. The red line is trajectory of penetrating gradient algorithm. After a few short movements, the line of gradient (accidentally) passed through a point close to the optimum but on the other side of the hill, far away from the current point. Gradient method did not see this point and could not reach it but penetrating gradient jumped right to it in a single step. The red line hits the hill at one side, penetrates through it and gets out on the other side, clearly demonstrating the penetrating property of the algorithm.

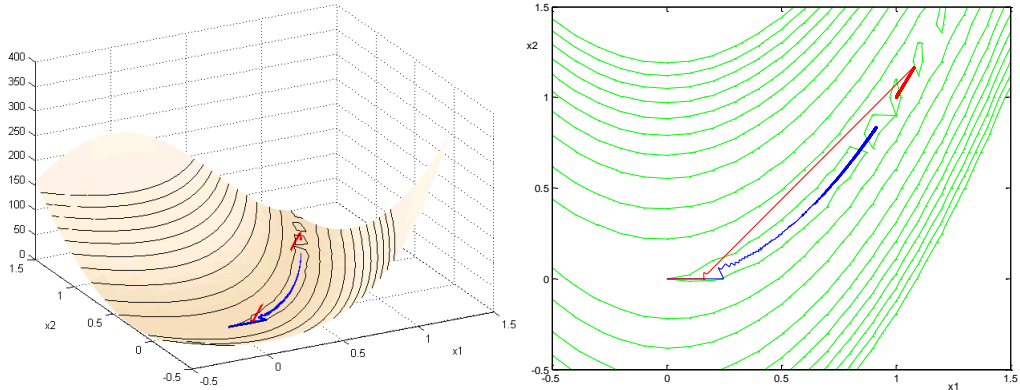

**Figure 1:** Comparison of gradient method that performs line search by backtracking and penetrating gradient algorithm that always searches in the direction of gradient. Blue line is trajectory of gradient method and red one of penetrating gradient algorithm. Left: surface of Rosenbrock function with levels (black lines) and algorithms' trajectories. Right: level sets (projections of levels on the $x_1x_2$ plane) with algorithms' trajectories.

**Benchmark tests**

The most interesting application of penetrating gradient is within the deepest descent method when it is used to examine several candidate search directions. In simulations that follow we compare three variants of Newton's method (with Levenberg-Marquardt modification) and *rss*1*rmax*2 algorithm. The three variants of Newton's method differ in the way of determination, i.e., adjustment of the step. They are: one that uses backtracking (backNLM), one using line search that satisfies strong Wolfe conditions (wolfeNLM) [1] and one applying penetrating gradient computation (Algorithm 1) but always only in Newton's direction (pgNLM). The deepest descent method is implemented as *rss*1*rmax*2 variant of Algorithm 2.

The selected algorithms were applied to many benchmark systems (taken from [22]) that have real solutions. The relevant data about this experiment is in Table 1 whose columns contain:

| | |
|---:|:---|
| name: | name of the system in corresponding reference |
| $n$: | dimension; number of equations, that is, unknowns |
| degree: | total degree of the system (product of degrees of all equations in the system) |
| starting point: | vector with initial values of all unknowns for which *rss*1*rmax*2 algorithm found a solution (zeros = all components initially zero, ones = all components initially one). Since the primary aim of this work is to examine the penetrating gradient algorithm, i.e., the deepest descent method, and different methods perform differently for the same starting point, we report only points from which *rss*1*rmax*2 algorithm found a solution, regardless of the results with other methods. It is mostly the case that only *rss*1*rmax*2 can solve more complicated systems in few attempts (starting from several points), though sometimes other methods can solve them starting from some other point. These points are, however, not specified in the table for conciseness. |
| iterations: | number of iterations performed by each algorithm to find a solution. Hyphen '-' denotes failure, i.e., the algorithm stopped due to some alternative stopping condition (for example, finding a local optimum), without solving the system. |

As can be seen in the table, $rss1rmax2$ algorithm solved all test systems (we consider system solved when $rss$ drops below $10^{-8}$). Occasionally this required several trials to find appropriate starting point but afterwards the solution was found very quickly. There is no execution time data in the table because it is machine dependent and it is proportional to the number of iterations, so number of iterations is the most relevant measure of performance. Of course, duration of a single iteration differs among the algorithms and depends on complexity of the system but penetrating gradient iteration is certainly the longest. Hence, to get a feeling about execution times it is sufficient to note that solving, for instance, *rabmo* system by $rss1rmax2$ algorithm on our computer[1] lasted about 20 seconds. Most of the systems were solved in a few iterations, which means practically instantaneously.

Several more examples can be found in [23] and they are all easy to solve by $rss1rmax2$ algorithm.

**Table 1:** Results for benchmark systems.

| name | *n* | degree | starting point | iterations $rss1rmax2$ | iterations backNLM | iterations wolfeNLM | iterations pgNLM |
|---|---|---|---|---|---|---|---|
| boon | 6 | 1024 | [0 0 1 1 0 0] | 8 | - | - | - |
| butcher | 7 | 4608 | [–1 –1 –1 0 –1 –1 –1] | 493 | - | 146 | 60 |
| butcher8 | 8 | 4608 | [1 2 3 4 5 6 7 8] | 176 | - | - | - |
| camera1s | 6 | 64 | [0 0 0 0 0 1] | 845 | - | - | - |
| caprasse | 4 | 144 | [–2 0 –2 2] | 4 | 6 | 7 | 4 |
| cassou | 4 | 1344 | [–50 0 0 0] | 367 | - | - | - |
| chemequs | 5 | 108 | zeros | 59 | 7 | 7 | 7 |
| cohn2 | 4 | 900 | ones | 4 | 15 | 15 | 3 |
| comb3000 | 10 | 96 | ones | 7 | 7 | 7 | 7 |
| cpdm5 | 5 | 243 | ones | 1 | 6 | 6 | 1 |
| cyclic7 | 7 | 5040 | zeros | 1 | 1 | 1 | 1 |
| cyclic8 | 8 | 40320 | [–2 –1 1 –1 –2 –2 1 2] | 24 | - | - | - |
| d1 | 12 | 4068 | [0 0 0 0 0 0 0 0 0 0 0 1] | 9 | 15 | 15 | 14 |
| des18_3 | 8 | 324 | [5 5 5 –5 –5 5 5 5] | 372 | - | - | - |
| des22_24 | 10 | 256 | zeros | 18 | 11 | 18 | 12 |
| discret3 | 8 | 256 | zeros | 5 | 6 | 6 | 5 |
| eco8 | 8 | 1458 | zeros | 2 | 2 | 2 | 2 |
| fourbar | 4 | 256 | ones | 13 | 39 | 37 | 90 |
| geneig | 6 | 243 | zeros | 8 | 9 | 7 | 8 |
| heart | 8 | 576 | [1 2 2 –3 –3 2 2 1] | 9 | 13 | 14 | 12 |
| i1 | 10 | 59049 | zeros | 2 | 2 | 2 | 2 |
| ipp | 8 | 256 | ones | 7 | - | - | - |
| katsura5 | 6 | 32 | zeros | 1 | 5 | 5 | 7 |
| kin1 | 12 | 4608 | ones | 4 | 5 | 5 | 4 |
| kinema | 9 | 64 | [1 2 0 –1 –1 –1 0 2 1] | 24 | 12 | 20 | 15 |
| ku10 | 10 | 1024 | zeros | 30 | - | - | - |
| lorentz | 4 | 16 | [1 2 2 1] | 4 | 4 | 4 | 4 |
| noon5 | 5 | 243 | zeros | 1 | 4 | 9 | 1 |
| proddeco | 4 | 256 | ones | 1 | 11 | 11 | 1 |
| puma | 8 | 128 | ones | 43 | 8 | 9 | 9 |
| quadfor2 | 4 | 24 | ones | 5 | 27 | 29 | 157 |
| quadgrid | 5 | 120 | [–5 –1 3 –1 –4] | 7 | - | - | 4 |
| rabmo | 9 | 36000 | [–1 .5 0 0 –1 0 0 .5 –1] | 1902 | - | - | - |
| rbpl | 6 | 486 | zeros | 20 | - | - | - |
| rbpl24es | 9 | 576 | zeros | 10 | 13 | 18 | 63 |
| redcyc6 | 5 | 120 | [2 1 1 –1 1 2] | 8 | 7 | - | 375 |
| redcyc7 | 6 | 720 | [0 0 0 1 0 0 0] | 1 | 1 | 1 | 1 |
| redcyc8 | 7 | 5040 | [1 0 1 1 1 0 0 1] | 935 | - | - | - |
| redeco6 | 6 | 16 | zeros | 4 | 4 | 4 | 4 |

---
[1] Intel Pentium Dual Core processor, ASUS P5K PRO motherboard, 2 GB RAM and Windows XP operating system.

| name | $n$ | degree | starting point | iterations $rss1rmax2$ | iterations backNLM | iterations wolfeNLM | iterations pgNLM |
|---|---|---|---|---|---|---|---|
| redeco7 | 7 | 32 | zeros | 4 | 4 | 4 | 4 |
| redeco8 | 8 | 64 | zeros | 4 | 4 | 4 | 4 |
| rediff3 | 3 | 8 | zeros | 4 | 5 | 5 | 4 |
| reimer5 | 5 | 720 | [–1 –1 1 1 –1] | 140 | - | - | - |
| rose | 3 | 216 | zeros | 2 | - | - | - |
| s9_1 | 8 | 16 | zeros | 11 | 9 | 11 | 23 |
| sendra | 2 | 49 | ones | 3 | 3 | 3 | 3 |
| solotarev | 4 | 36 | ones | 1 | 1 | 1 | 1 |
| trinks | 6 | 24 | zeros | 15 | 5 | 5 | 7 |
| virasoro | 8 | 256 | ones | 1 | 5 | 5 | 1 |
| wood | 4 | 36 | zeros | 4 | 12 | 82 | 104 |
| wright | 5 | 32 | zeros | 1 | 3 | 4 | 1 |

## 5. Penetrating gradient algorithm for quasi-linear and linear systems

As a subclass of polynomial systems, defined by the constraint that maximal power of variables in the system (2.7) is one, i.e., the terms are of the form $a \cdot x_1 x_2 \ldots x_n$, quasi-linear systems lend themselves to penetrating gradient algorithm straightforwardly and enable important numerical and program simplifications. Search in arbitrary directions can be achieved in the same way as with full polynomial systems, meaning that we have to construct polynomial (2.11) and find all its roots to identify the proper one. But, if search direction is an axis direction, that is, if $\mathbf{d}$ is a unit vector, then all arguments of $P_i(\mathbf{x}+\alpha \cdot \mathbf{d})$ in (2.8) remain unchanged (simply $x_i$), except the one corresponding to the selected axis, being $(x_j+\alpha)$. This is easily seen in example (2.9) by substituting, for instance, $d_1 = 1$ and $d_2 = 0$, which would mean that the search direction is $x_1$ axis. Thus, after substitution of $x_i+\alpha \cdot d_i$ into terms like $x_1 x_2 \ldots x_n$, they remain almost the same, just instead of the selected variable, let us say $x_1$, they will contain $(x_1+\alpha)$ and be something like $(x_1+\alpha) \cdot x_2 \ldots x_n$. Since all $x_2 \ldots x_n$ are known, the new terms of the form $(x_1+\alpha) \cdot x_2 \ldots x_n$ will be simply $(x_1+\alpha)$ multiplied by a number. Because, in this example, the search direction is $x_1$ axis, we search for new (better) $x_1$, that is, we need $(x_1+\alpha)$. However, it is equal whether we fix $x_1$ and calculate $\alpha$ or we simply calculate $(x_1+\alpha)$, i.e., new $x_1$, at once. This leads to especially simple form of penetrating gradient algorithm presented in this section.

Let us denote the selected unknown (i.e, coordinate axis corresponding to the search direction) as $x_j$ and introduce $\lambda = x_j + \alpha$ for what we want to compute. After substitution of known variables into the system (2.7), the only unknown will be $\lambda$ and we obtain a system like

$$\begin{aligned} P_1(\lambda) &= a_{11} \cdot \lambda + a_{12} \cdot \lambda + \ldots + a_{1q} \cdot \lambda = c_1 \\ P_2(\lambda) &= a_{21} \cdot \lambda + a_{22} \cdot \lambda + \ldots + a_{2q} \cdot \lambda = c_2 \\ &\ldots \\ P_m(\lambda) &= a_{m1} \cdot \lambda + a_{m2} \cdot \lambda + \ldots + a_{mq} \cdot \lambda = c_m \quad, \end{aligned} \qquad (5.1)$$

which is nothing but

$$\begin{aligned} (a_{11} + a_{12} + \ldots + a_{1q}) \cdot \lambda &= a_1 \cdot \lambda = c_1 \\ (a_{21} + a_{22} + \ldots + a_{2q}) \cdot \lambda &= a_2 \cdot \lambda = c_2 \\ &\ldots \\ (a_{m1} + a_{m2} + \ldots + a_{mq}) \cdot \lambda &= a_m \cdot \lambda = c_m \end{aligned}$$

or, with $\mathbf{a} = [a_1 \, a_2 \, \ldots \, a_m]^\text{T} \in \mathrm{R}^m$ and $\mathbf{c} = [c_1 \, c_2 \, \ldots \, c_m]^\text{T} \in \mathrm{R}^m$, simply

$$\mathbf{a} \cdot \lambda = \mathbf{c} \quad . \qquad (5.2)$$

System (5.1), i.e., (5.2) is a linear system of $m$ equations in one unknown $\lambda$ and its LS solution can be relatively simply computed as [17],[19],[20],

$$\lambda_{\text{LS}} = (\mathbf{a}^\text{T} \mathbf{a})^{-1} \mathbf{a}^\text{T} \mathbf{c} = \frac{\sum_{i=1}^{m} a_i c_i}{\sum_{i=1}^{m} a_i^2} = \frac{\sum_{i=1}^{m} a_i c_i}{\|\mathbf{a}\|^2} \quad . \qquad (5.3)$$

Thus, if search direction is an axis direction, we can find new (better) estimate of solution without construction of the derivative of summary polynomial (2.11) and computing its roots. This also means that there is no need for computing derivatives of system equations, which is a considerable simplification compared to the procedure required for solving polynomial systems. Potential problem that $\|\mathbf{a}\| = 0$ can easily be prevented by a simple test in the program.

Thereby we arrive to the simplest, derivative free, variant of penetrating gradient algorithm, suitable specifically for quasi-linear and linear systems.

**Algorithm 3:** Derivative-free penetrating gradient algorithm for quasi-linear and linear systems
1. Given the current point **x** and one of the coordinate axes as search direction (let it be $j^{th}$ axis), consider the corresponding variable $x_j$ unknown.
2. In system (2.7) substitute current values for all variables considered known and denote $x_j$ as $\lambda$. This will render the system into a linear system in $\lambda$ of the form (5.1), i.e., (5.2).
3. If $\|\mathbf{a}\| \neq 0$, solve (5.2) for $\lambda$ according to (5.3). Let us denote the result as $\lambda_{LS}$. If $\|\mathbf{a}\| = 0$, just return some indication of failure to the calling routine.
4. Compose new point $\mathbf{x}_{new}$ by replacing $x_j$ in **x** by $\lambda_{LS}$. This $\mathbf{x}_{new} = [x_1 \ldots x_{j-1}\ \lambda_{LS}\ x_{j+1} \ldots x_n]$ minimizes the *rss* of the original system (2.7) in selected dimension (i.e., along the line of selected coordinate axis).
5. Check stop conditions and if none is satisfied repeat the procedure from the step 1.

In [24] and [25] there is a more elaborate study of solving quasi-linear and linear systems using this algorithm and it shows to be very effective even in this simplest form. However, there is a surprising property of the algorithm when applied to pure linear systems. It turns out that, when applied to a pure linear system $\mathbf{A}\mathbf{x} = \mathbf{b}$ in a way that all axes are taken as search directions in a sequence, Algorithm 3 is nothing but another form of Gauss-Seidel method [16] for solving linear systems applied to the system of normal equations $\mathbf{A}^T\mathbf{A}\mathbf{x} = \mathbf{A}^T\mathbf{b}$. Thus, we can establish the following result.

**Theorem:** *When applied to a linear system* $\mathbf{A}\mathbf{x} = \mathbf{b}$, *where* $\mathbf{A} \in R^{n \times n}$ *and* $\mathbf{b} \in R^n$, *in a way that all axes are taken as search directions in a sequence, Algorithm 3 is equivalent to Gauss-Seidel method applied to the system* $\mathbf{A}^T\mathbf{A}\mathbf{x} = \mathbf{A}^T\mathbf{b}$. □

*Proof*: We prove the statement by showing that penetrating gradient computation of new value of a single unknown is algebraically equivalent to Gauss-Seidel computation.

Penetrating gradient, i.e., Algorithm 3 applied to $\mathbf{A}\mathbf{x} = \mathbf{b}$ tries to reduce $rss(\mathbf{x}) = \mathbf{r}(\mathbf{x})^T\mathbf{r}(\mathbf{x}) = \|\mathbf{b} - \mathbf{A}\mathbf{x}\|^2$ by changing only one $x_k$ at a time. To this aim, it fixes all other unknowns and solves $d/dx_k [rss(x_k)] = 0$. Expanded, the condition $d/dx_k [rss(x_k)] = 0$ is

$$2(\mathbf{b} - x_1 \cdot \mathbf{a}_1 - \ldots - x_k \cdot \mathbf{a}_k - \ldots - x_n \cdot \mathbf{a}_n)^T \cdot \mathbf{a}_k = 0, \quad (5.4)$$

where $\mathbf{a}_j$ means $j^{th}$ column of matrix **A**. By isolating $x_k$ in (5.4) we see that new $x_k$ by Algorithm 3 is

$$x_k = \frac{(\mathbf{b} - x_1 \mathbf{a}_1 - \ldots - x_{k-1}\mathbf{a}_{k-1} - x_{k+1}\mathbf{a}_{k+1} - \ldots - x_n \mathbf{a}_n)^T \mathbf{a}_k}{(\mathbf{a}_k)^T \mathbf{a}_k} \quad (5.5)$$

If we introduce notation $\mathbf{C} = \mathbf{A}^T\mathbf{A}$ and $\mathbf{d} = \mathbf{A}^T\mathbf{b}$ and have in mind that matrix **C** is symmetric, formula (5.5) can be written as

$$x_k = \frac{d_k - x_1 c_{k1} - \ldots - x_{k-1}c_{k,k-1} - x_{k+1}c_{k,k+1} - \ldots - x_n c_{kn}}{c_{kk}} \quad (5.6)$$

On the other hand, new value of $y_k$ yielded by Gauss-Seidel method applied to the system $\mathbf{C}\mathbf{y} = \mathbf{d}$ is given by

$$\sum_{i=1}^{k-1} y_i c_{ki} + y_k c_{kk} + \sum_{i=k+1}^{n} y_i c_{ki} = d_k \quad (5.7)$$

where from we get

$$y_k = \frac{d_k}{c_{kk}} - \sum_{i=1}^{k-1} y_i c_{ki} - \sum_{i=k+1}^{n} y_i c_{ki}$$

$$= \frac{d_k - y_1 c_{k1} - \ldots - y_{k-1}c_{k,k-1} - y_{k+1}c_{k,k+1} - \ldots - y_n c_{kn}}{c_{kk}} \quad (5.8)$$

Comparing formula (5.6) for new $x_k$ yielded by Algorithm 3 with formula (5.8) for new $y_k$ yielded by Gauss-Seidel method, we see that they are identical, so for any starting point $\mathbf{x}^{(0)} = \mathbf{y}^{(0)}$ the sequence of new values of unknowns produced by either method will be identical as well. ∎

The equivalence of Algorithm 3 and Gauss-Seidel method sheds a new light on Gauss-Seidel method but more importantly, it enables transferring all theoretical results developed for Gauss-Seidel method to Algorithm 3. An important Algorithm 3 property follows as the corollary of the previous theorem.

**Corollary:** *When applied to a linear system $\mathbf{Ax} = \mathbf{b}$, where $\mathbf{A} \in \mathrm{R}^{n \times n}$ and $\mathbf{b} \in \mathrm{R}^n$, in a way that all axes are taken as search directions in a sequence, Algorithm 3 converges to a least squares solution $\mathbf{x} = (\mathbf{A}^T\mathbf{A})^{-1}\mathbf{A}^T\mathbf{b}$ of the system.* □

*Proof*: *By Theorem, Algorithm 3 is equivalent to Gauss-Seidel method applied to the system $\mathbf{A}^T\mathbf{Ax} = \mathbf{A}^T\mathbf{b}$. It is well known from linear algebra that matrix $\mathbf{A}^T\mathbf{A}$ is positive definite and since the convergence of Gauss-Seidel method to a least squares solution when applied to positive definite systems is proven, it follows that Algorithm 3 converges to a least squares solution as well. Least squares solutions of a system $\mathbf{Ax} = \mathbf{b}$ are $\mathbf{x} = (\mathbf{A}^T\mathbf{A})^{-1}\mathbf{A}^T\mathbf{b}$, which proves the corollary.* ∎

Unfortunately, precise analysis of computational complexity [24],[25] reveals that penetrating gradient iteration requires exactly three times more basic operations (additions, multiplications) than Gauss-Seidel iteration, so there is no much reason to use Algorithm 3 for solving linear systems in practice. It is true that Algorithm 3 can be used as "a universal" LS solver of linear systems, while Gauss-Seidel method imposes certain constraints on the system structure, but if we wish LS solution of a linear system there are other methods (like QR or SVD, see [15],[16]) at disposal and, with full systems, they will likely be faster than Algorithm 3. With sparse systems Algorithm 3 might be competitive but this strongly depends on implementation and is subject to a comprehensive examination out of the scope of this paper.

## 6. Complexity of penetrating gradient algorithm

Obviously, the complexity of penetrating gradient algorithm depends on two major steps: construction of polynomial (2.11) and computation of its roots. Complexities of both these steps depend on the degrees $\beta_i$ of equations in the system (see discussion below (2.11)) and these degrees depend on the structure of particular terms (powers of unknowns in the first place). Therefore, to make any analysis feasible, let us assume that the highest degree in the whole system is given in advance and denote it $\beta$. Then, in the worst case, the construction of (2.11) requires $m$ multiplications of $\beta$-degree and $(\beta-1)$-degree polynomials and summation of the results. For simplicity, we shall neglect summation since it is much faster than multiplication. If polynomial multiplication is performed as time domain convolution [29], which is recommended with respect to the processing speed as long as polynomial degrees do not exceed few tenths, one multiplication will be $O(\beta+(\beta-1)-1)$ operation so the construction of the whole polynomial (2.11) will require $O(m \cdot [\beta+(\beta-1)-1]) = O(m \cdot 2\beta-2) = O(m \cdot \beta)$ operations (for "big Oh" notation see [18] or [2]). Computation of roots of polynomial (2.11) is much more demanding and is numerically one of the most challenging computational tasks. The details are not the subject of this paper and interested reader should consult the literature on numerical analysis, for example [14]. In this approximate analysis we can take that computation of roots of polynomial is roughly $O(\beta^2)$ operation [27], so we can assert that the overall complexity of penetrating gradient is

$$O(m \cdot \beta + \beta^2). \tag{6.1}$$

This is per direction cost, so when penetrating gradient is used to check several directions in one iteration as in *rss1rmax2* algorithm, the cost of iteration will grow proportionally to the number of directions searched. Because complexity of penetrating gradient is not directly related to the parameters of the objective function and its gradient, precise comparison with theoretical complexities of other methods is not liable. In any case, it is obvious that penetrating gradient algorithm is numerically the most demanding one but it is also true that the benefits it provides are well worth of increased computational effort.

## 7. Conclusion

The deepest descent method powered by penetrating gradient algorithm efficiently solves polynomial systems of up to moderate complexity. The penetrating gradient algorithm replaces line search and increases probability of discovering a global optimum, i.e., solution of the system due to the two distinguishing properties: the penetrating property and the ability to search in arbitrary directions. The penetrating property, that is, the fact that the algorithm yields the best point on the whole line of search direction is itself a remarkable advantage over any line search routine because line search can find new point only in some limited area around the current point. However, the full power of the algorithm emanates from the mutual aid of penetration and search along arbitrary directions and this is what provides the opportunity for new strategies like the deepest descent. The deepest descent strategy is not suitable for common line search routines because they are not effective with non-descent search directions in the current point. In contrast, for penetrating gradient algorithm search direction is unimportant, so even directions that are not descent directions looking from the current point are appropriate and can be beneficial. Thus, the determination of the best point on the line of search direction is now solved and all efforts can be directed to the recognition of the best search directions. It can be expected that this will soon yield many new ideas and strategies.

Results of simulations show that the deepest descent method is more successful than classical optimization methods and the cost paid for this improvement is relatively low. It is true that penetrating gradient requires more computation than other algorithms but it is very suitable for parallelization because one of its most time-consuming steps is construction of (2.11), that is, the derivative of summary polynomial. Polynomial (2.11) is nothing but a sum of $m$ independent polynomials and each of them can be processed in a separate program thread. Thus, properly parallelized, construction of polynomial (2.11) can take just about $2\beta$ operations.

A difficulty with penetrating gradient algorithm is the need for determination of the roots of a relatively high-degree polynomial. Moreover, at present we compute all the roots, although we actually need only one of them. The problem is that we still cannot determine in advance which root is "the right one" and this is certainly one of the challenges for the further research. Solving this issue would significantly reduce computational requirements of penetrating gradient and possibly improve its numerical stability, which would, in turn, enable solving more complex systems than at present.

Another problem is that at present we cannot control the algorithm's progress and force it to find all solutions to the system. Although it would probably be possible to restrict the search to a limited area by imposing constraints on the acceptable values of unknowns, we still do not have a reliable manner to direct the algorithm so that it avoids already searched paths. What is needed is some kind of path separation control like in homotopy continuation method and this is certainly the second main subject of the further development.

Finally, the presented computation procedure works only with real numbers and the algorithm can find only real solutions to the system. This can be as much advantage as disadvantage depending on the application, but there is no theoretical obstacle for extending the computation to the complex domain so this is certainly an important aim in the near future.